\documentclass[english,a4paper,twoside,12pt]{amsart}

\usepackage{graphicx}

\usepackage[english]{babel}
\usepackage[mathscr]{eucal}
\usepackage[latin1]{inputenc}

\def\phi{\varphi}
\def\Bbb{\mathbb}
\def\Cal{\mathscr}

\renewcommand{\frak}[1]{\mathfrak{#1}}

\theoremstyle{plain}
\newtheorem{thm2}{Theorem}[section]
\newtheorem{lem2}[thm2]{Lemma}
\newtheorem{prop2}[thm2]{Proposition}
\newtheorem{cor2}[thm2]{Corollary}

\theoremstyle{definition}
\newtheorem{ex2}[thm2]{Example}
\newtheorem{rem2}[thm2]{Remark}

\newcommand{\lem}[1]{\begin{lem2} #1 \end{lem2}}

\DeclareMathOperator{\Id}{Id}

\DeclareMathOperator{\ad}{ad}

\newcommand{\Js}{J_{\frak s}}
\newcommand{\Jg}{J_{\frak g}}

\begin{document}

\title[Levi--Tanaka extensions] {Abelian extensions of \\
semisimple graded $CR$ algebras}

\author[A. Altomani]{Andrea Altomani}
\address{Scuola Normale Superiore \\ Piazza dei Cavalieri, 7 \\ 
56100 Pisa (Italy)}
\email{altomani@sns.it}
\curraddr{Dipartimento di Matematica \\ Universit\`a di Roma ``Tor Vergata''\\
Viale della Ricerca Scientifica snc\\ 00133 Roma (Italy)}

\author[M. Nacinovich]{Mauro Nacinovich}
\address{Dipartimento di Matematica \\ Universit\`a di Roma ``Tor Vergata''\\
Viale della Ricerca Scientifica snc\\ 00133 Roma (Italy)}
\email{nacinovi@mat.uniroma2.it}

\date{February 24th, 2003}

\begin{abstract}
In this paper we take up the problem of describing the $CR$ vector
bundles $M$ over compact standard $CR$ manifolds $S$, which 
are themselves standard
$CR$ manifolds. They are associated to 
special graded Abelian
extensions of semisimple graded $CR$ algebras.
\end{abstract}



\maketitle

\tableofcontents

\section{Introduction}
\par
The Levi--Tanaka algebras were introduced by Noburu Tanaka in
\cite{Tan67,Tan70} to study  the differential geometry
of $CR$ manifolds. These
algebras were 
further investigated in \cite{MN97}, where also
the special class of {\it standard} $CR$ manifolds associated to them
was introduced. In the paper \cite{MNssempl} all semisimple
Levi--Tanaka algebras were classified, and later it was shown that
they correspond to special minimal orbits for actions of real
Lie groups on flag manifolds (see \cite{MN2000,MNVigo,Wf69}),
which in turn are precisely the
standard $CR$ manifolds that are compact. For nonstandard homogeneous
$CR$ manifolds the problem of classification is still open. Recently
that of classifying the 
compact homogeneous $CR$ manifolds of hypersurface
type ($CR$ codimension $1$) has been solved in \cite{AS02}.
\par
In this paper we take up the problem of describing the $CR$ vector
bundles $M$ over compact standard $CR$ manifolds $S$, which 
are themselves standard
$CR$ manifolds. They are naturally associated to graded Abelian
extensions of semisimple graded $CR$ algebras. These
extensions are defined in \S 2; in \S 3 we compute the
Lie algebra $\frak g$ of infinitesimal $CR$ automorphisms of $M$; 
then, after
showing that the partial complex structure of $M$ is defined
by an inner derivation of $\frak g$, we obtain necessary and
sufficient conditions in order that a graded Abelian extension of
a semisimple graded $CR$ algebra $\frak s$ admits a $CR$ structure.
The last section is devoted to the discussion of some examples. 

\section{Abelian extensions of graded Lie algebras}

\subsection{Representations of graded Lie algebras}
\par
Let  $\frak s=\oplus_{p\in\Bbb Z}{\frak s_p}$ be a finite dimensional 
$\Bbb Z$-graded real Lie algebra. Given  
a finite dimensional $\frak s$-module $\frak l$, 
the Abelian extension of $\frak s$ by $\frak l$ is
$\frak s\oplus\frak l$, with Lie product 
defined by:
\[
[X,Y]=\begin{cases} [X,Y]_{\frak 
s}&\quad\text{if}\quad X,Y\in\frak s\, ,\\
X\cdot Y &\quad\text{if}\quad X\in\frak s\, ,\; Y\in\frak l\\
-Y\cdot X&\quad\text{if}\quad Y\in\frak s\, ,\; X\in\frak l\\
0&\quad\text{if}\quad X,Y\in\frak l\, .
\end{cases}
\]
We will use 
interchangeably the words {\it module} and {\it representation}.
\par\smallskip
Given $Y\in\frak l$ and $X_1,...,X_r\in \frak s$ we define by
recurrence:
\[
\begin{cases}
[Y]=Y\, &\quad\text{if}\quad r=0\, ,\\
[X_1,Y]=X_1\cdot Y&\quad\text{if}\quad r=1\, ,\\
[X_1,X_2,Y]=[X_1,[X_2,Y]]&\quad\text{if}\quad r=2\, ,\\
[X_1,X_2,...,X_r,Y]=[X_1,[X_2,...,X_r,Y]]
&\quad\text{if}\quad r>2\, .
\end{cases}
\]
We call $[X_1,X_2,...,X_r,Y]$ a \emph{Lie monomial} in $Y$ of length $r$;
we say that  
it is \emph{homogeneous} if $X_1,...,X_r$ are homogeneous elements of 
$\frak s$; homogeneous \emph{decreasing} (resp.  \emph{increasing}) 
if moreover 
$\deg(X_1)\geq\deg(X_2)\geq\cdots\geq\deg(X_r)$ (resp.  
$\deg(X_1)\leq\deg(X_2)\leq\cdots\leq\deg(X_r)$); the integer 
$\deg(X_1)+\cdots +\deg(X_r)$ 
is called its {\it degree}.  The monomial $[Y]$  
is homogeneous of degree $0$.
\par

\lem{\label{lem:liemonom}Let $\frak s$ be a graded real Lie algebra,
$\frak l$
an $\frak s$-module, and $Y\in\frak l$. \par The 
homogeneous decreasing 
(or increasing) Lie monomials in $Y$ generate
the $\frak s$-submodule $\frak l'$ of $\frak l$ generated by $Y$.}

\begin{proof} Clearly $\frak l'$ is generated by
the homogeneous Lie monomials in $Y$. 
Thus, to prove our contention, it suffices to show that every
Lie monomial in $Y$, homogeneous of degree $d$, is a linear 
combination of
homogeneous decreasing (or increasing) Lie monomials in $Y$
of the same degree $d$. 
Let $r>1$, let $[X_1,...,X_r,Y]$ be a homogeneous Lie monomial and
$(i_1,...,i_r)$ a permutation of the indices
 $(1,...,r)$.
We claim that $[X_1,...,X_r,Y]-[X_{i_1},...,X_{i_r},Y]$ is
a linear combination of homogeneous Lie monomials in $Y$ of the 
same
degree $d$ and of
length $<r$.
Since the group of permutations 
of $\{1,2,...,r\}$ is generated by the transpositions exchanging
 $i$ and $(i+1)$ for $1\leq i\leq r-1$, this follows from the formula:
\begin{multline*} 
[X_1,\hdots,X_i,X_{i+1},\hdots,X_r,Y]=\\
=[X_1,\hdots,X_{i+1},X_i,\hdots,X_r,Y]+[X_1,\hdots,[X_i,X_{i+1}],
\hdots,X_r,Y]\, ,
\end{multline*}
where we note that the homogeneous Lie monomials in the right hand 
side
have the same degree of the homogeneous Lie monomial in the left
hand side.
This proves our claim. The statement of the
Lemma follows then by recurrence
on the length $r$ of a homogeneous Lie monomial $[X_1,...,X_r,Y]$.
\par
\end{proof}
\par\medskip
A \emph{gradation} of the $\frak s$-module $\frak l$ is a decomposition of
 $\frak l$ into a direct sum of finite dimensional vector subspaces:
 $\frak l=\oplus_{p\in\Bbb Z}{\frak l_p}$ such that
\[
[\frak s_p,\frak l_q]\subset\frak l_{p+q}\qquad
\forall p,q\in\Bbb Z\, .
\]
\par
\medskip
If $\frak a$ is a graded vector space, we will denote by $\frak a_{-}$ 
the subspace $\bigoplus_{p<0}\frak a_{p}$ and by $\frak a_{+}$ 
the subspace $\bigoplus_{p\geq 0}\frak a_{p}$ (note that $\frak 
a_{0}\subset\frak a_{+}$ but $\frak 
a_{0}\cap\frak a_{-}=\{0\}$).  If $\frak a$ is a 
$\Bbb Z$-graded Lie algebra, then 
$\frak a_{+}$ and $\frak a_{-}$ are Lie subalgebras of $\frak a$.
\par
\smallskip
A graded representation $\frak l$ of the graded Lie algebra
$\frak s$ is called 
\begin{itemize}
\item \emph{transitive} if $[\frak s_{-1},Y]\neq \{0\}$ for $Y\in\frak
l_p\setminus \{0\}$ and $p\geq 0$;
\item \emph{nondegenerate} if $[\frak s_{-1},Y]\neq \{0\}$ for $Y\in\frak
l_{-1}\setminus \{0\}$; 
\item \emph{fundamental} if $\frak l_{p-1}=[\frak 
s_{-1},\frak l_p]$ for $p<0$, i.e. if $\frak l_{-1}$ 
generates $\frak l_-$ as an $\frak s_-$-module.
\end{itemize}
Note that if $\frak l\neq\{0\}$ 
is transitive and fundamental, then $\frak 
l_{-1}$ is different from $\{0\}$.
\par
\smallskip
A graded Lie algebra $\frak s$ is called
\emph{transitive} (resp.\ \emph{nondegenerate}, \emph{fundamental}) if 
it is 
\emph{transitive} (resp.\ \emph{nondegenerate}, 
\emph{fundamental}) 
when considered as a graded $\frak s$-module via the adjoint representation.
\par
\smallskip
If $\frak s=\oplus_{p\in\Bbb Z}{\frak s_p}$ is a graded Lie algebra,
the linear map $D:\frak s \rightarrow\frak s$ defined by
$$D(X)=pX\quad {\rm if}\quad p\in\Bbb Z\quad {\rm and}\quad X\in\frak s_p$$
is a derivation of $\frak s$. When this derivation $D$ is inner,
we say that $\frak s$ is \emph{characteristic} and call
characteristic any element $E_{\frak s}$ of $\frak s$ such that
$D={\rm ad}_{\frak s}(E_{\frak s})$.
\begin{lem2}
Every semisimple graded Lie algebra $\frak s$ is characteristic
and contains a unique characteristic element $E_{\frak s}$.
\end{lem2}
\begin{proof}
Indeed, all derivations of a semisimple Lie algebra are inner.
Moreover the characteristic element $E_{\frak s}$ is unique
because the adjoint representation is faithful.
\end{proof}
\par
Recall that
the \emph{kind} and \emph{co-kind} of a finite dimensional $\Bbb Z$-graded Lie
algebra $\frak s=\oplus_{p\in\Bbb Z}{\frak s_p}$ are the integers
$\mu=\mu(\frak s)=\sup\{p\in\Bbb Z \mid  \frak s_{-p}\neq 0\}$ and
$\nu=\nu(\frak s)=\sup\{p\in\Bbb Z\mid\frak s_{p}\neq 0\}$.
If $\frak s$ is semisimple, then $\mu(\frak s)=\nu(\frak s)$.
\par
Likewise the kind $\mu(\frak l)$ and the co-kind $\nu(\frak l)$ of
a finite dimensional graded $\frak s$-module 
$\frak l=\oplus_{p\in\Bbb Z}{\frak l_p}$ are defined by
$\mu(\frak l)=\sup\{p\in\Bbb Z \mid  \frak l_{-p}\neq 0\}$ and
$\nu(\frak l)=\sup\{p\in\Bbb Z\mid\frak l_{p}\neq 0\}$.
\par
\begin{lem2}
Let $\frak s$ be a graded Lie algebra, and $\frak l$ a graded 
representation of
$\frak s$. Then
\begin{itemize}
\item[($i$)] $\frak{s_-\oplus l_-}$ is fundamental if and only if 
$\frak s_-$
and $\frak l_-$ are both fundamental;
\item[($ii$)] if $\frak s$ and $\frak l$ are both transitive (resp.\ 
nondegenerate) then $\frak{s\oplus l}$ is transitive (resp.\ 
nondegenerate);
\item[($iii$)] if $\frak{s\oplus l}$ is transitive (resp.\ 
nondegenerate) then
$\frak l$ is transitive (resp.\ nondegenerate);
\end{itemize}
\end{lem2}
\begin{proof}
The proof is straightforward.
\end{proof}

\lem{\label{lem:grad}
Assume that $\frak s$ is a semisimple 
graded real Lie algebra.
Then every finite dimensional $\frak s$-module $\frak l$ admits
a gradation. If moreover $\frak l$ is irreducible, and
 $\frak l=\oplus_{p\in\Bbb Z}{\frak l_p}=\oplus_{p\in\Bbb Z}{\frak 
l'_p}$
are two gradations of $\frak l$, then there exists an integer $k$
such that $\frak l'_p=\frak l_{p+k}$.}

\begin{proof}
 Assume that $\frak l$ is irreducible.
Let $E_{\frak s}$ be the characteristic element of 
 $\frak s=\oplus_{-\mu\leq p\leq\mu}{\frak s_p}$.
Then $T:\frak l\ni Y \to [E_{\frak s},Y]\in\frak l$ 
is a semisimple endomorphism
with rational eigenvalues. The difference of any two eigenvalues of
 $T$ is an integer. Hence, if we fix an eigenvalue $\lambda$ of $T$
and define 
$\frak l_h=\{Y\in\frak l \mid  [E_{\frak s},Y]=(\lambda+h)Y\}$
for all $h\in\Bbb Z$ we obtain a gradation 
 $\frak l=\oplus_{h\in\Bbb Z}{\frak l_h}$ of the $\frak s$-module 
$\frak l$.

Vice versa, if $\frak l=\oplus_{h\in\Bbb Z}{\frak l_h}$ is a
gradation of the irreducible $\frak s$-module $\frak l$,
the subspaces $\frak l_h$ are $T$-invariant because
$E_{\frak s}\in\frak s_0$. 
Hence they are generated by eigenvectors of $T$. Suppose that 
 $\frak l_d$, for some $d\in\Bbb Z$, contains a nonzero vector
 $Y$. We can assume that $Y$ is an eigenvector of $T$. If
 $[E_{\frak s},Y]=\lambda Y$, 
we have $\lambda\in\mathbb Q$ and, by choosing
a basis of the real vector space $\frak l$ consisting of
homogeneous Lie polynomials in $Y$, we obtain that
  $\frak l_h=\{Z\in\frak l \mid [E_{\frak s},Z]=(\lambda+h-d)Z\}$.
\par
In general, after decomposing $\frak l$ into a direct sum of
irreducible $\frak s$-submodules, we obtain a gradation of $\frak l$ 
by
defining a gradation on each irreducible $\frak s$-submodule
of the decomposition.\par
This completes the proof of the lemma. 
\end{proof}
\par
\begin{prop2}
\label{prop:gradeddecomp}
Let $\frak s$ be a graded semisimple Lie
algebra. Every graded $\frak s$-module $\frak l$ admits 
a decomposition into a direct sum of graded irreducible
$\frak s$-modules.
This decomposition is unique up to $0$-degree isomorphisms.
\par
A graded representation of $\frak s$ 
is transitive (resp.\ nondegenerate, 
fundamental) if and only if all its irreducible graded components are 
transitive (resp.\ nondegenerate, fundamental).
\end{prop2}

\begin{proof}
Let $D:\frak s\oplus\frak l\rightarrow\frak s\oplus\frak l$
be the derivation of $\frak s\oplus\frak l$ defined by
$D(X)=pX$ if $X\in\frak s_p\oplus\frak l_p$,
for $p\in\Bbb Z$. 
The derivation $D$ commutes with
the action of the characteristic element $E_{\frak s}$ of $\frak s$.
We identify $\frak s$, via the adjoint representation, to
a subalgebra of the Lie algebra $\frak D$ of derivations
of $\frak s\oplus\frak l$ and denote by
$\frak s'$ the Lie subalgebra of $\frak D$ generated by
$\frak s$ and $D$. Since $D-E_{\frak s}$ is $0$ on $\frak s$,
$D-E_{\frak s}$ and $\frak s$ commute. Thus $\frak s'$ is reductive
with center $\Bbb R\cdot (D-E_{\frak s})$.\par
Next we note that $\frak l$ is $\frak s'$-invariant and
$D-E_{\frak s}$ is diagonalizable on $\frak l$.
Thus the representation 
of $\frak s'$ on $\frak l$ is semisimple, and 
$\frak l$  decomposes
into a direct sum of $D$-invariant, i.e. graded, irreducible
$\frak s$-modules. This decomposition is unique, modulo isomorphisms
of $\frak s'$-modules. These isomorphisms, commuting with $D$,
are of degree zero.
\par
The last statement is straightforward.
\end{proof}
\par

\subsection{Graded $CR$ algebras and representations}

A \emph{graded $CR$ 
algebra} is a graded real Lie algebra 
$\frak s=\oplus_{p\in\Bbb Z}{\frak s_p}$,
with a complex structure $J:\frak s_{-1}\to\frak s_{-1}$
satisfying: 
\begin{itemize}
\item[($i$)] $J^2=-\Id_{\frak s_{-1}}$; 
\item[($ii$)] $[JX,JY]=[X,Y]$ for every $X$ and $Y$ 
in $\frak s_{-1}$;
\item[($iii$)] $[A,JX]=J[A,X]$ for every $A\in\frak s_0$ and $X\in\frak s_{-1}$.
\end{itemize}
A graded representation $\frak l$ of $\frak s$ is called a 
\emph{$CR$ representation} if a complex structure is given on 
$\frak l_{-1}$, that will also be
denoted by $J$, making the Abelian extension 
$\frak s\oplus\frak l$ a graded $CR$ algebra.
\par

\begin{lem2}\label{lem:unique}
Given a graded $CR$ Lie algebra $\frak s$ and a nondegenerate
graded representation $\frak l$, there is at most one $CR$ structure
on $\frak l$ that makes $\frak{s\oplus l}$ into a $CR$ extension of $\frak
s$.
\end{lem2}
\begin{proof}
Let $J_{\frak l}$ and $J'_{\frak l}$ be two complex structure on 
$\frak
l_{-1}$ that make $\frak{s\oplus l}$ into a $CR$ extension of $\frak
s$. Then, if $Y\in\frak l_{-1}$, we 
have 
\[
[J_{\frak l}Y,X]=-[Y,J_{\frak s}X], \qquad 
[J'_{\frak l}Y,X]=-[Y,J_{\frak s}X]
\]
for every $X\in\frak s_{-1}$, thus
\[
[J_{\frak l}Y-J'_{\frak l}Y,X]=0
\]
for every $X\in\frak s_{-1}$. Because of the non degeneracy of $\frak 
l$ we
obtain that $J_{\frak l}Y=J'_{\frak l}Y$ for every $Y\in\frak l_{-1}$.
\end{proof}
\par
\begin{prop2}\label{prop:jlinl}
Let $\frak l$ be a transitive nondegenerate graded $CR$ representation of 
a semisimple graded $CR$ algebra $\frak s$.  Then every graded
$\frak s$-submodule 
of $\frak l$ is also $CR$.
\end{prop2}

\begin{proof}
Let $\frak l'$ be a graded $\frak s$-submodule of $\frak l$, $Y\in\frak 
l'_{-1}$ a nonzero vector and $\frak l''$ a graded $\frak s$-submodule 
complementary to $\frak l'$.  Write $JY = Y'+Y''$ with 
$Y'\in\frak l'$ and $Y''\in\frak l''$. 
For every $X\in\frak s_{-1}$, we have
$[X,Y']+[X,Y'']=[X,JY]=-[JX,Y]\in\frak l'$.  
This implies that $[X,Y'']=0$ for every $X\in\frak s_{-1}$. 
By the assumption 
that $\frak l$ is non degenerate, we obtain that $Y''=0$, and 
therefore $JY=Y'\in\frak l'$.
\end{proof}

\section{Levi--Tanaka extensions}
\par
We recall the definition of a Levi--Tanaka algebra. If 
$\frak m=\bigoplus_{p<0}\frak m_{p}$ 
is a fundamental nondegenerate
graded $CR$ algebra then there exists a unique 
(modulo isomorphisms) transitive graded finite 
dimensional $CR$ algebra 
$\frak g$, maximal with respect to inclusion, such that 
$\frak m=\frak g_{-}$.
This maximal prolongation
$\frak g=\frak g(\frak m)$ is called a Levi--Tanaka 
algebra. 
\par
In fact we can take
$\frak g=\oplus_{p\in\Bbb Z}{\frak g_p}$ with the
subspaces $\frak g_p$ recursively defined by:
\[
\frak g_p=\begin{cases}
{\frak m_p}&\quad {\rm if}\quad p<0\, ;\\
\{A\in\frak D_0(\frak m,\frak m)\, | \, A(JX)=JA(X)\;
\forall X\in\frak m_{-1}\}\, &\quad {\rm if}\quad p=0\, ;\\
\frak D_p(\frak m,\oplus_{q<p}{\frak g_q})
&\quad {\rm if}\quad p>0\, ;
\end{cases}
\] 
where $\frak D_p$ denotes
degree $p$ homogeneous derivations.
\par
We refer to \cite{Tan67, Tan70, MN97, MN2000} for a thorough 
discussion of 
Levi--Tanaka algebras.  \par 
It is worth rehearsing
that a transitive semisimple 
graded {$CR$} algebra $\frak s$ is a Levi--Tanaka algebra if and
only if 
is nondegenerate.  This
is equivalent to the fact that every simple ideal of $\frak s$ has 
kind $\geq 2$.  
\par \smallskip
Let $\frak s$ be a semisimple graded {$CR$} 
algebra and let $\frak l$ be an $\frak s$-module.  Assume that $\frak 
l$ is graded, transitive and {$CR$}.  The Abelian extension $\frak 
s\oplus\frak l$ has a natural structure of graded {$CR$} algebra.  
Under the further assumption that $\frak s\oplus\frak l$ is 
fundamental and nondegenerate, we consider the maximal transitive 
{$CR$} prolongation $\frak g=\frak g(\frak s\oplus\frak l)$ of $\frak 
g_{-}=\frak s_-\oplus\frak l_-$.  It 
is a Levi--Tanaka algebra containing $\frak s\oplus\frak l$.  \par We 
will say that $\frak{s}\oplus\frak{l}$ is the \emph{partial 
Levi--Tanaka 
extension} of $\frak{s}$ by $\frak l$ 
and $\frak g$ is the \emph{Levi--Tanaka 
extension} of $\frak s$ by $\frak l$.
\par
\subsection{The structure of the Levi--Tanaka prolongation}
Let $\frak s$ 
be a semisimple graded {$CR$} algebra, and 
$\frak l$ a graded, 
transitive, {$CR$}, fundamental and 
nondegenerate $\frak s$-module.
Let $\frak g=\oplus_{p\in\Bbb Z}{\frak g_p}$ be the Levi--Tanaka
extension of $\frak s$ by $\frak l$.
Denote by
\[
\pi:\frak s\oplus\frak l\to\frak l
\]
the projection 
onto $\frak l$ along $\frak s$.
\par
\begin{lem2}
$\pi$ is a $0$-degree derivation of $\frak s\oplus\frak l$, commuting
with $\frak s_0$.
\end{lem2}
\par
Since $\pi(\frak g_{p})=\frak l_p\subset\frak g_p$ for all $p<0$ 
and $\pi$ commutes with the partial complex structure $J$ on $\frak 
g_{-1}=\frak s_{-1}\oplus\frak l_{-1}$, there is a unique element in 
$\frak g_0$, that we still denote by $\pi$, such that $[\pi,X]=\pi(X)$ 
for all $X\in\frak s\oplus\frak l$.  \par \lem{$\ad_{\frak 
g}(\pi):\frak g \to \frak g$ is semisimple.}
\par
\begin{proof} 
Indeed $\ad_{\frak g}(\pi)|_{\frak g_{-}}$ is 
semisimple and 
hence the statement follows from Lemma 3.8 in \cite{MN97}.
\end{proof}
\par
We say that a partial Levi--Tanaka extension $\frak s\oplus\frak 
l$ is \emph{semisimple} if the corresponding Levi--Tanaka extension 
$\frak g(\frak s\oplus\frak l)$ is semisimple, and that $\frak 
s\oplus\frak l$ and $\frak g(\frak s\oplus\frak l)$ are \emph{proper} 
if $\frak g(\frak s\oplus\frak l)$ does not contain any semisimple 
ideal.
\par
\begin{prop2}\label{prop:properdec}
Let $\frak s$ be a semisimple fundamental graded {$CR$} algebra and 
$\frak l$ a 
transitive fundamental graded {$CR$} $\frak s$-module, such that 
$\frak s\oplus\frak l$ is nondegenerate.  Let $\frak g$ be the 
Levi--Tanaka extension of $\frak s$ by $\frak l$.\par 
Then there exist two graded 
{$CR$} ideals $\frak s'$ and $\frak s''$ of $\frak s$ and two graded 
{$CR$} 
$\frak s$-submodules $\frak l'$ and $\frak l''$ of $\frak l$ such 
that:
\begin{itemize}
\item[($i$)] $\frak s=\frak s'\oplus\frak s''$, $ \frak l=\frak 
l'\oplus 
\frak l''$;
\item[($ii$)] $\frak s'\oplus\frak l'$ and $\frak s''\oplus\frak l''$ 
are ideals in $\frak s\oplus\frak l$;
\item[($iii$)] $\frak g=\frak g(\frak s'\oplus\frak l')\oplus\frak 
g(\frak 
s''\oplus\frak l'')$ as a sum of ideals;
\item[($iv$)] $\frak g(\frak s'\oplus\frak l')$ is a semisimple 
Levi--Tanaka extension of $\frak s'$;
\item[($v$)] $\frak g(\frak s''\oplus\frak l'')$ is a proper 
Levi--Tanaka 
extension of $\frak s''$.
\end{itemize}
\end{prop2}
\par
\begin{proof}
Let $\frak r$ 
be the radical of $\frak g$.
Since $\frak r\cap\frak s=0$, the map
$\ad(\pi)|_{\frak r_{-}}:\frak r_{-}\to \frak l_{-}$ 
is injective. This implies that $\frak r_{-}$ 
is contained in $\frak l_{-}$. Let $\frak l''=\frak l\cap\frak r$, so 
that $\frak l''_{-}=\frak r_{-}$.
\par
Let $\frak L$ be a graded {$CR$} Levi subalgebra of $\frak g$, 
containing $\frak s$ and
${\rm ad}_{\frak g}(\pi)$-invariant (its existence is granted by 
\cite{MNVigo}).  Let $\frak l'=\frak L\cap\frak l$.  Then $\frak 
l=\frak l'\oplus\frak l''$, because $\frak l_{-}=\frak 
l'_{-}\oplus\frak l''_{-}$ and every irreducible component of $\frak 
l$ is generated by its degree $-1$ homogeneous elements.  Every 
simple 
ideal $\frak a$ of $\frak s$ is contained in a simple ideal $\frak 
L_{\frak a}$ of $\frak L$.  Define
\[
\frak s'=\bigoplus\left\{\frak a\mid \frak a \text{ is 
a simple ideal of $\frak s$ and }\frak L_{\frak a}\cap\frak l\neq 
0\right\}\oplus\ker_{\frak s}\frak l
\]
and
\[
\frak s''=\bigoplus\left\{\frak a\mid \frak a \text{ is 
a simple ideal of $\frak s$, $\frak L_{\frak a}\cap\left(\frak l
\oplus\ker_{\frak s}\frak l\right)= 0$} \right\}.
\]
Then condition ($i$) is fulfilled and condition ($ii$) 
follows because 
$[\frak s',\frak l'']=0$ and $[\frak s'',\frak l']=0$;
condition ($iii$) is a consequence of ($ii$) and Proposition~3.3 in 
\cite{MN97}; 
condition ($iv$) holds because $\frak g(\frak s'\oplus\frak 
l')=\bigoplus_{\frak a\subset\frak s'}\frak L_{\frak a}$ is a direct 
sum of simple ideals.  \par Finally, assume by contradiction that 
$\frak g(\frak s''\oplus\frak l'')$ contains a a simple ideal $\frak 
a$.  Then $[\frak a,\frak l'']=0$ because $\frak l''\subset\frak r$, 
and $[\frak a,\frak l']=0$ because $\frak a$ is not contained in 
$\frak s'$.  Hence $\frak a\subset \ker_{\frak s}\frak l\subset\frak 
s'$ and we get a contradiction.  The proof is complete.
\end{proof}
\par
\begin{ex2}\label{ex:sl2c1}
Let
\[
\frak s=\frak{sl}(2,\Bbb C)=\left.\left\{ \left(\begin{matrix} 
a_{11}&a_{12}\\
a_{21}&a_{22}
\end{matrix}\right)\right|\, a_{ij}\in\Bbb C,\; 
a_{11}+a_{22}=0\right\}.
\]
We consider the gradation $\frak s=\oplus_{-1\leq p\leq 1}{\frak s_p}$
corresponding to the characteristic element 
\[
E_{\frak s}=
\left(\begin{matrix}
1/2&0\\
0&-1/2
\end{matrix}\right)
\]
and the {$CR$} structure defined by the matrix
\[
J_{\frak s}=
\left(\begin{matrix}
i/2&0\\
0&-i/2
\end{matrix}\right).
\]
\par We consider the space $\frak l$ of anti-Hermitian $2\times 2$ 
matrices and we let $\frak s$ act on $\frak l$ by:
\[
[X,A]=X\cdot A\, =\, AX+X^*A\qquad\text{for $X\in\frak s$, $A\in\frak 
l$}\, .
\]
We define on $\frak l$ the gradation corresponding to the 
characteristic
element 
\[
[E,A]=[E_{\frak s},A]-A\quad \text{for $A\in\frak l$}\, 
\]
and the {$CR$} structure obtained by restriction to $\frak l_{-1}$
of the action of $J_{\frak s}$ on $\frak l$.\par
Define $\frak s\frak u(2,2)$ as the simple Lie algebra of
$4\times 4$ complex matrices $Y$ satisfying
\[
YI_{2,2}+I_{2,2}Y^*=0
\]
for
\[
I_{2,2}=\left(
\begin{matrix}
0&0&1&0\\
0&0&0&1\\
1&0&0&0\\
0&1&0&0
\end{matrix}\right)\, .
\]
Then $\frak s\oplus\frak l$ is isomorphic to the Lie subalgebra of
$\frak s\frak u(2,2)$ consisting of matrices of the form
\[
\left(\begin{matrix}
X&A\\
0&-X^*
\end{matrix}\right)\qquad
\text{with $X\in\frak s$ and $A\in\frak l$.}
\]
and clearly $\frak g(\frak s\oplus\frak l)=\frak s\frak u(2,2)$ is
semisimple.
\end{ex2}
\par
\begin{ex2}
Let $\frak s=\frak{sl}(2,\Bbb C)$ with the same gradation and $CR$ 
structure
as in the previous example and let $\frak l=\frak l'\oplus\frak l''$ 
where
both $\frak l'$ and $\frak l''$ are isomorphic to the representation 
$\frak l$
of the previous example, with the same gradation and the same $CR$ 
structure. 
\par
In this case $\frak g(\frak s\oplus\frak l)=\frak s\oplus\frak 
l\oplus\Bbb R\pi'\oplus\Bbb R\pi''$, where $\pi'$ and $\pi''$ are 
projections onto $\frak l'$ and $\frak l''$ respectively.
In this case $\frak g$ is proper.
It is the Levi--Tanaka algebra of 
kind two and $CR$ codimension two,
associated to a singular pencil of $3\times 3$ Hermitian matrices
(see \cite{MN2001}).\par
\end{ex2}
\par
\medskip
Next we give a sufficient condition for a Levi--Tanaka
extension to
be proper.
\par
\begin{lem2}\label{lemsb}
Let $\frak g$ be a finite dimensional Levi--Tanaka 
algebra and let $\frak s$ be a semisimple Levi--Tanaka subalgebra of 
$\frak g$.  Let $\phi_{-}:\frak g_{-}\to \frak s_{-}$ be a 
homomorphism of graded {$CR$} algebras, whose restriction to $\frak 
s_{-}$ is the identity.  Then $\phi_{-}$ extends to a split 
surjective 
homomorphism of graded {$CR$} algebras
\[
\tilde\phi:\frak g \to \frak s\, .
\]
\end{lem2}
\par
\begin{proof} 
For $A\in\frak g_0$, the map 
$\phi_0(A):\frak s_-\to\frak s_-$ defined by:
\[
\phi_0(A)(X)=\phi_{-}([A,X])\qquad\forall X\in\frak s_-\, 
\]
is a degree zero derivation of $\frak s_{-}$,
that commutes with the complex structure
$J_{\frak s}$ on $\frak s_{-1}$,  and hence defines
an element of $\frak s_0\subset\frak D_0(\frak s_-,\frak s_-)$.
This gives a homomorphism $\phi_0:\frak g_0\to\frak s_0$.
\par
Next we define recursively
$\phi_p:\frak g_p\to\frak s_p$, for $p>0$,  by
\[
\phi_p(A)(X)=\phi_{p-h}([A,X])\qquad
\forall A\in\frak g_p,\quad\forall X\in\frak s_{-h},\quad h>0\, .
\]
One verifies that 
$\phi(A)$ is an element of 
$\frak s_p=\frak D_p(\frak s_{-},\oplus_{q<p}{\frak s_q})$.
\par
Define $\tilde\phi:\frak g\to\frak s$ 
to be equal to $\phi_{-}$ on $\frak g_{-}$ and to
$\phi_p$ on $\frak g_p$ for $p\geq 0$. Since $\tilde\phi$
is the identity on $\frak s\subset\frak g$, the map
$\tilde\phi:\frak g\to \frak s$ is onto and splits. \par
It is clear that $\tilde\phi$ is a surjective homomorphism
of graded {$CR$} algebras.
\end{proof}
\par
\begin{prop2}
If $\frak s$ is a semisimple Levi--Tanaka algebra and  $\frak l$ a 
nondegenerate
faithful graded {$CR$} $\frak s$-module,
then the Levi--Tanaka extension $\frak g$
of $\frak s$ by $\frak l$ is proper.
\end{prop2}
\par
\begin{proof}
From the previous lemma we know that there is a graded {$CR$} ideal 
$\frak b$ of 
$\frak g$ such that $\frak g=\frak s\oplus\frak b$.
Since $[\pi,\frak b]\subset\frak b\cap\frak l$, we obtain that
$\frak l_-$, and hence $\frak l$, are contained in $\frak b$.
Assume by contradiction that $\frak g$ contains a simple ideal
$\frak a$.  
\par
If $[\pi,\frak a]\neq\{0\}$ or $\frak a\cap\frak b\neq\{0\}$ then 
$\frak 
a\subset\frak b$.  Therefore $\frak a_{-1}\subset\frak l_{-1}$ and 
$[\frak a_{-1},\frak s_{-1}\oplus\frak l_{-1}] \subset [\frak a,\frak 
s]+[\frak l,\frak l]=\{0\}$.  By the condition that $\frak 
s_-\oplus\frak l_-$ is nondegenerate, we obtain that $\frak a_{-1}$, 
and hence $\frak a_-$, are $\{0\}$.  This contradicts the 
transitivity 
of $\frak g$.
\par
If $\frak a\cap\frak b=\{0\}$ and $[\pi,\frak a]=0$ then $\frak 
a_{-}\subset\frak s_{-}$, but this contradicts the faithfulness of 
$\frak l$, because $[\frak a,\frak b]=\{0\}$.
\end{proof}
\par
If $\frak g$ is a proper Levi--Tanaka extension of $\frak s$
by a transitive fundamental nondegenerate 
graded {$CR$} $\frak s$-module $\frak l$,  
Proposition \ref{prop:properdec} and Lemma \ref{lemsb} yield a decomposition: 
\[
\frak g=\frak s\oplus\frak b=\frak s\oplus\frak a\oplus\frak r
\]
where $\frak b$ is an ideal containing $\frak l$, 
which is the direct sum of the
radical $\frak r$ of $\frak g$ and of a semisimple 
$\frak a\subset\frak g_0$ that is an ideal in the Levi
subalgebra $\frak L=\frak s\oplus\frak a$ of $\frak g$.
\par
\medskip
Denote by $T:\frak b \to \frak b$ the restriction of 
$\ad_{\frak g}(\pi)$ to $\frak b$. Since $T$ is
semisimple, its minimal polynomial $m_T(\lambda)$
is a product $m_T(\lambda)=p_1(\lambda)\cdots p_t(\lambda)$
of irreducible distinct polynomials in $\Bbb R[\lambda]$.
Set:
\[
\frak b^{p_i}=\{X\in\frak b\, | \, p_i(T)(X)=0\}\qquad
\text{for}\quad i=1,...,t\, .
\]
Being $\frak l\subset\frak b$, we can assume 
that $\frak l\subset\frak b^{p_1}$, with 
$p_1(\lambda)=\lambda-1$. 
Moreover, all $\frak b^{p_i}$ are graded, with kind
$\mu(\frak b^{p_i})\leq 0$ if $i\geq 2$, because
$\frak b_-=\frak l_-\subset\frak b^{p_1}$.
The $\frak b^{p_i}$'s are $\frak s$-modules, as $T$ commutes
with the restriction of $\ad_{\frak g}(\frak s)$ to $\frak b$.
For every $2\leq i\leq t$, the $\frak s$-module $\frak b^{p_i}$
contains an element $U$ homogeneous of minimal degree $d\geq 0$.
Hence $[U,\frak s_{-1}]=\{0\}$ and there exists $Y\in\frak l_{-1}$
such that $[U,Y]\neq 0$. This implies that for each $i$ with 
$2\leq i\leq t$, also
the polynomial 
$p_i(\lambda-1)$ is a factor of the minimal polynomial $m_T(\lambda)$
of $T$, because $p_i(T-I)([U,Y])=0$. From this
observation  we deduce that all 
irreducible factors of $m_T(\lambda)$ are of the first order and
 $T$ can be diagonalized:
its eigenvalues are integers $h$ with $2-t\leq h\leq 1$.
Thus we rewrite:
\[
\frak b=\bigoplus_{2-t\leq h\leq 1}{\frak b^h}
\]
where 
\[
\frak b^h=\{X\in\frak b\, | \, [\pi,X]=hX\,\},\qquad h\in\Bbb Z,\;
2-t\leq h\leq 1\, 
\]
are the eigenspaces of $T$.\par
Moreover, $\frak b^1=\frak l$. Indeed, if this were not the case,
there would be a nonzero
graded $\frak s$-submodule $\frak w$ of $\frak b^1$ which is a 
complement
of $\frak l$ in $\frak b^1$. Take a nonzero
element $W$ homogeneous of minimal degree in $\frak w$. Then,
because $\frak g$ is transitive, there is $X\in\frak l_{-1}$
with $[W,X]\neq 0$. But this gives a contradiction, because
$T([W,X])=2[W,X]$, and all eigenvalues of $T$ are $\leq 1$.
\par\smallskip
The ideal $\frak b$ of $\frak g$ contains the 
radical $\frak r$ of $\frak g$.  In particular the radical 
$\frak r(\frak b)$ of $\frak b$ coincides with the radical $\frak r$ 
of $\frak g$. Thus the decomposition 
$$\frak b=\frak a\oplus\frak r$$
is a Levi--Malcev decomposition of $\frak b$. Since $\pi$ is a
semisimple element of $\frak b_0$, we can choose
$\frak a\subset\frak b_0$ to be $\pi$-invariant. Thus we
decompose $\frak a$ into a direct sum of eigenspaces of $T$:
$\frak a=\oplus_{h\leq 0}{\frak a^h}$.
Since $\oplus_{h<0}{\frak a^h}$ is a nilpotent ideal of $\frak a$,
it is $\{0\}$ and we obtain 
$\frak a=\frak a^0$, i.e. $[\pi,\frak a]=0$.
\par
\smallskip
We have obtained that
\[
[\frak a,\frak s]=\{0\}\quad\text{and}\quad [\frak a,\frak 
l]\subset\frak l\, .
\]
Therefore $\frak a$ is a semisimple
subalgebra of the algebra of zero degree
$\frak s$-endomorphisms of $\frak l$. 
The discussion beolow will show that actually $\frak a$ is
a Levi subalgebra of the algebra of zero degree
$\frak s$-endomorphisms of $\frak l$. 
\par
\smallskip
Now observe that $\frak r$
(considered as a subalgebra of $\frak g\frak l_{\Bbb R}(\frak g)$ via
the adjoint representation) is splittable and therefore decomposes 
into
the direct sum of the ideal $\frak n(\frak g)$ consisting of its
nilpotent elements and of a maximal Abelian
subalgebra $\frak t$ of $\frak r(\frak g)$ consisting of semisimple
elements. Since $[\pi,\frak L]=\{0\}$, the semisimple element $\pi$ 
belongs to $\frak r$. Thus we can take
$\pi\in\frak t$ and $\frak t\subset\frak g_0$. \par
The nilpotent ideal $\frak n(\frak g)$ 
is $\ad_{\frak g}(\pi)$-invariant and therefore
decomposes into eigenspaces of $T$:
\[
\frak n(\frak g)=\bigoplus_{\substack{h\in\Bbb Z\\h\leq 0}}{\frak 
n^h(\frak g)}
\] 
where $\frak n^h(\frak g)=\{X\in\frak n(\frak g)\, | \,
 [\pi,X]=hX\}$.
Moreover $\frak n^1(\frak g)=\frak l$ and $\frak n^0(\frak g)$ is a
subalgebra and an $\frak s$-submodule, with 
 $[\frak n^0(\frak g),\frak l]\subset\frak l$. 
\par
\bigskip
\paragraph{\textbf{The irreducible case}}
Now we suppose that $\frak l$ is irreducible. In this case we know 
from Ch.~9 of \cite{Bourbaki} that for $\frak a$ we have the cases:
\begin{itemize}
\item[($i$)] the representation $\frak l$ is of the real or of the
complex type: then $\frak a=\{0\}$;
\item[($ii$)] the representation $\frak l$ is of the quaternionic 
type:
then $\frak a\simeq\frak s\frak u(2)$.
\end{itemize}
By Engel's theorem the set
of $Y\in\frak l$ such that $[\frak n^0(\frak g),Y]=0$ is different
from $\{0\}$. Since it is an $\frak s$-submodule of $\frak l$, and we 
are
supposing that $\frak l$ is irreducible, it
must coincide with $\frak l$. From this we derive that 
$[\frak n^0(\frak g),\frak l]=\{0\}$. But $\frak n^0(\frak g)$ is
graded and $\frak n^0(\frak g)=\bigoplus_{p\geq 0}{\frak n^0_p(\frak 
g)}$.
For a nonzero homogeneous term $U$ of minimal degree in $\frak 
n^0(\frak g)$
we would have $[U,\frak s_{-1}]=\{0\}$ and $[U,\frak l_{-1}]=0$,
contradicting the transitivity of $\frak g$. Thus
$\frak n^0(\frak g)=\{0\}$ and this implies that $\frak n^h(\frak 
g)=\{0\}$
for all $h\leq 0$.

Hence we obtain $\frak n(\frak g)=\frak l$ and $\frak t=\frak 
r^0(\frak g)$
consists of semisimple elements. We have $\frak t=\Bbb R\pi$ if 
$\frak l$
is of the real or of the quaternionic type and  $\frak t=\Bbb C\pi$ if
$\frak l$ is of the complex type.
\par
We have proved the following:
\begin{thm2}\label{thm:irredstruct}
Let $\frak s$ be a semisimple graded {$CR$} algebra.  Let $\frak l$ 
be an 
irreducible  nondegenerate {$CR$} graded $\frak s$-module.
Assume that the Levi--Tanaka extension $\frak g$
of $\frak s$ by $\frak l$ is proper. 
Then $\frak g$ is a finite 
dimensional 
real Lie algebra and we have:
\begin{itemize}
\item[(\emph{i})] if $\frak l$ is of the real type, 
then $\frak s$ is a Levi subalgebra of $\frak g$, and
$\frak g$ admits the {$CR$}
Levi-Malcev decomposition:
\[
\frak g=\frak s\oplus\left(\frak l\oplus\Bbb R\pi\right),
\]
where $\frak l$
is the maximal nilpotent ideal of the radical $\frak r$ of $\frak g$ 
and $\Bbb R\pi$ a maximal Abelian subalgebra of semisimple elements
of $\frak r$;
\item[(\emph{ii})] if $\frak l$ is of the complex type, 
then $\frak s$ is a Levi subalgebra of $\frak g$, and
$\frak g$ admits the {$CR$}
Levi-Malcev decomposition:
\[
\frak g=\frak s\oplus\left(\frak l\oplus\Bbb C\pi\right),
\]
where $\frak l$
is the maximal nilpotent ideal of the radical $\frak r$ of $\frak g$ 
and $\Bbb C\pi$ a maximal Abelian subalgebra of semisimple elements
of $\frak r$;
\item[(\emph{iii})] If $\frak l$ is of the quaternionic type,
$\frak g$ has a
{$CR$}
Levi-Malcev decomposition:
\[
\frak g=\left(\frak s\oplus\frak a\right)\oplus
\left(\frak l\oplus\Bbb R\pi\right),
\]
with a Levi subalgebra $\frak L=\frak s\oplus\frak a$ consisting
of the direct sum of $\frak s$ and of a
simple ideal $\frak a\subset\frak g_0$, 
isomorphic to $\frak s\frak u(2)$;  $\frak l$
is the maximal nilpotent ideal of the radical $\frak r$ of $\frak g$,
and $\Bbb R\pi$ is a maximal Abelian subalgebra of semisimple elements
of $\frak r$.
\end{itemize}
\end{thm2}
\par
\bigskip
\paragraph{\textbf{The general case}}
If the representation $\frak l$ is not irreducible, the structure of 
$\frak g$
may be more complicated (see the example below). 
\par
\smallskip
Decompose $\frak l$ into a direct sum of irreducible {$CR$} graded 
representations:
\[
\frak l = \bigoplus_{i} m_i\cdot\frak l^i
\]
where the $m_i$'s are positive integers and 
the $\frak l^i$'s irreducible 
graded
{$CR$} representations of $\frak s$, pairwise not isomorphic as 
graded 
representations (i.e. distinct $\frak l^i$'s 
can be isomorphic, but endowed with
non isomorphic gradations).
\par
We keep the notation used throughout this section. In particular,
$\frak a\subset\frak g_0$ is the Levi part 
of the ideal $\frak b$ complementary to $\frak s$ in $\frak g$ and 
$\frak t\subset\frak g_0$ 
is a maximal Abelian subalgebra of semisimple elements of 
the radical $\frak r$ of $\frak g$.\par
Then 
$$\frak{a\oplus t}=\bigoplus_i(\frak a_i\oplus\frak t_i)$$ 
where 
$(\frak a_i\oplus\frak t_i)$ 
is isomorphic to $\frak{gl}(m_i,\mathbb K_i)$ and
$\mathbb K_i$ is equal to $\mathbb R$, $\mathbb C$ or $\mathbb H$ 
in agreement with
the type of the representation $\frak l^i$. The action of $(\frak
a_i\oplus\frak t_i)$ on $m_i\cdot\frak l^i$ is induced by the 
canonical action
of $\frak{gl}(m_i,\mathbb K_i)$ on the span of the maximal weight 
vectors of
each copy of $\frak l^i$. Thus we have the cases:
\begin{itemize}
\item[($i$)] $\frak l^i$ is of the real or the complex type: then 
$\frak
a_i\simeq\frak{sl}(m_i,\mathbb K_i)$ and $\frak t^i\simeq\mathbb K_i$;
\item[($ii$)] $\frak l^i$ is of the quaternionic type: then $\frak
a_i\simeq\frak{sl}(m_i,\mathbb H_i)\oplus\frak{su}(2)$ and $\frak
t^i\simeq\mathbb R$. 
\end{itemize}
\par
We summarize the discussion above in the following:
\begin{thm2}
Let $\frak s$ be a transitive semisimple graded {$CR$} algebra and 
$\frak 
l$
a nondegenerate {$CR$} graded $\frak s$-module. 
\par
Assume that the Levi--Tanaka extension $\frak g$ of $\frak s$
by $\frak l$ is proper.\par
 Let $\frak l=\bigoplus m_i\cdot\frak l^i$ be a 
decomposition of $\frak l$ into irreducible graded $\frak s$-modules.
Denote by 
$\pi_i$ the projection onto $m_i\cdot\frak l^i$.  Then $\frak g$ 
is finite 
dimensional and admits the {$CR$} Levi-Malcev decomposition:
\[
\frak g = (\frak s \oplus \bigoplus_i \frak a^i) \oplus
(\frak l \oplus \bigoplus_i \frak t^i \oplus \frak n),
\] 
where, for every $i$, $\frak a^i$ and $\frak t^i$ are contained in
$\frak g_0$ and we distinguish the following cases:
\begin{itemize}
\item[($i$)] $\frak l^i$ is of the real type: 
then $\frak a_i\simeq\frak{sl}(m_i,\mathbb R)$ and $\frak t^i=\mathbb 
R\pi_i$;
\item[($ii$)] $\frak l^i$ is of the complex type: 
then $\frak a_i\simeq\frak{sl}(m_i,\mathbb C)$ and $\frak t^i=\mathbb 
C\pi_i$;
\item[($iii$)] $\frak l^i$ is of the quaternionic type: 
then $\frak a_i\simeq\frak{sl}(m_i,\mathbb H)\oplus\frak{su}(2)$ and 
$\frak t^i=\mathbb R\pi_i$.
\end{itemize}
\end{thm2}
\par
Unlike the irreducible case, if $\frak l$ is reducible
the nilpotent ideal $\frak n$ 
can be different from $\{0\}$:\par
\begin{ex2}\label{example}
Let $\frak s=\frak{sl}(3,\mathbb C)$. It admits an essentially  unique
structure of a Levi--Tanaka algebra (see~\cite{MNssempl}), with
gradation and partial complex structure defined respectively by
the matrices
\[
E_{\frak s}=\begin{pmatrix}
-1&0&0\\
0&0&0\\
0&0&1
\end{pmatrix}
\qquad\text{and}\qquad
\Js=\begin{pmatrix}
\frac{1}{3}i&0&0\\
0&-\frac{2}{3}i&0\\
0&0&\frac{1}{3}i
\end{pmatrix}\, .
\]
\par
Then $\frak s_-$ is the subalgebra of nilpotent upper
triangular matrices, and $\frak s_0$ is the Cartan subalgebra of 
diagonal matrices of $\frak s$. 
\par
Let $\frak l$ be the direct sum of the standard representation $\frak 
l'=\mathbb C^3$
with basis 
$(e_i)_{1\leq i\leq 3}$ and its dual $\frak l''=(\mathbb C^3)^*$ with
dual basis 
$(f_j)_{1\leq j\leq 3}$. Choose the gradations of $\frak l'$ and 
$\frak l''$
by imposing that $\deg e_1=-3$, $\deg e_2=-2$, $\deg e_3=-1$, $\deg 
f_1=0$, $\deg f_2=-1$ and $\deg f_3=-2$. 
Then $\frak{s_-\oplus l_-}$ is a fundamental 
nondegenerate
graded Lie algebra. Multiplication by $i=\sqrt{-1}$ 
defines the {$CR$} structure on $\frak l$. Denote by
$\frak{g=g(s_-\oplus l_-)}$ the corresponding Levi--Tanaka 
algebra. 
In this case $\frak l$ is an ideal, and via the adjoint 
representation we can identify $\frak{s\oplus t\oplus n}$ to a 
subalgebra of $\frak{gl(l)}\simeq\frak{gl}(6,\Bbb C)$.
\par
Under this identification we have:
\begin{align*}
\frak s=&\left.\left\{
\begin{pmatrix}
A&0\\
0&-{^t\!A}
\end{pmatrix}\,\right|\, A\in\frak{sl}(3,\Bbb C)\right\}\, ,\\
\frak t=&\left.\left\{
\begin{pmatrix}
\alpha I&0\\
0&\beta I
\end{pmatrix}\,\right|\,
I=\begin{pmatrix}1&0\\0&1\end{pmatrix},\,\alpha,\beta\in\Bbb 
C\right\}\, ,\\
\intertext{and}
\frak n=&\left.\left\{
\begin{pmatrix}
0&0\\
N&0
\end{pmatrix}\,\right|\, N\in\frak{o}(3,\Bbb C)\right\}\, , 
\end{align*}
graded in such a way that
\[
N_0=
\begin{pmatrix}
0 & 0 & 0 \\
0 & 0 &-1 \\
0 & 1 & 0
\end{pmatrix},\quad
N_1=
\begin{pmatrix}
0 & 0 & -1 \\
0 & 0 & 0 \\
1 & 0 & 0
\end{pmatrix}\quad\text{and}\quad
N_2=
\begin{pmatrix}
0 & -1 & 0 \\
1 &  0 & 0 \\
0 &  0 & 0
\end{pmatrix}
\]
correspond to elements of degrees $0$, $1$ and $2$, respectively.
\end{ex2}

\subsection{The $CR$ structure 
of the Levi--Tanaka extensions}\label{sect:CRstruct}
\par
Our next goal is to show that the Levi--Tanaka algebras $\frak g$
that are prolongations of nondegenerate graded Abelian $CR$ extensions
of semisimple $CR$ algebras have the $J$ property: this means that
the partial complex structure $J$ on $\frak g_{-1}$ is the restriction
to $\frak g_{-1}$ of the inner derivation with respect to an element 
$J_{\frak g}$ of $\frak g$.
\par
All semisimple graded $CR$ algebras have the $J$ property
(see \cite{MNssempl}). Thus we can restrict our consideration to the
case where $\frak g$ is proper.\par
If $\frak s=\oplus_{p\in\Bbb Z}{\frak s_p}$ is a 
fundamental transitive semisimple
graded $CR$ algebra,
we denote by $J_{\frak s}$ the element of $\frak s_0$
such that $[J_{\frak s},X]=JX$ for all $X\in\frak s_{-1}$.
\par
\bigskip
We shall prove the following:
\par
\begin{thm2}\label{teo:Jprop}
Let $\frak s$ be a transitive semisimple graded $CR$ algebra.
Let $\frak l$ be a nondegenerate transitive 
$CR$ graded $\frak s$-module. Then the Levi--Tanaka extension 
$\frak g$ of $\frak s$ by $\frak l$ has the $J$ property.
\end{thm2}
\par
\medskip
The proof will be divided in several steps.
\par
\medskip

\begin{lem2}\label{lem:multidentity}
Let $\frak s$ be a transitive semisimple graded $CR$ algebra and 
$\frak l$
an irreducible nondegenerate graded $CR$ $\frak s$-module.
Assume that 
the Levi--Tanaka extension $\frak g$ of $\frak s$
by $\frak l$ is proper. Then $\frak g$  has the $J$ property if
and only if 
$\ad_{\frak g}(J_{\frak s})$ acts as a multiple of
the identity on $\frak l_{-2}$.
\par
Moreover if the $\frak s$-module $\frak l$ is real or quaternionic then
$J_{\frak g}=\Js$;  if $\frak l$ is an $\frak s$-module of the
complex type then $\Jg=\Js+ik\pi$, with $k\in\Bbb R$.
\end{lem2}
\par
\begin{proof}
Assume that there is an element $J_{\frak g}$
of $\frak g$ such that
$JX=[J_{\frak g},X]$ for $X\in\frak g_{-1}
=\frak s_{-1}\oplus\frak l_{-1}$.
Since $\ad_{\frak g}(J_{\frak g})$ 
is semisimple (see \cite{MN97}) and agrees with
$\ad_{\frak g}(J_{\frak s})$ on 
$\frak s_{-1}$, then $\frak s$ is $J_{\frak g}$-invariant and
we have $[J_{\frak g}-J_{\frak s},\frak s]=0$.
Hence,
$Z=J_{\frak g}-J_{\frak s}$ defines an 
$\frak s$-endomorphism of $\frak l$ which, by Schur's lemma,
is a multiple of the identity:
\[
\ad_{\frak g}(Z)= k \pi, \quad k\in\Bbb K,
\]
where $\Bbb K$ are 
either the real or the complex numbers, or the quaternions,
according to the fact that $\frak l$ is real, complex or
quaternionic. 
\par
Assume vice versa that $[J_{\frak s},X]=kX$, for some $k\in\Bbb K$,
when $X\in\frak l_{-2}$, and define: 
\[
J_{\frak g}=J_{\frak s}-k\pi\, .
\]
Then we have, for $X\in\frak l_{-1}$ and $Y\in\frak s_{-1}$:
\[
[[J_{\frak g},X],Y]=[J_{\frak g},[X,Y]]+[X,[J_{\frak g},Y]]=
[X,[J_{\frak s},Y]]=[X,JY]\, .
\]
This shows that $J_{\frak g}$ defines the complex structure of
$\frak g$, so that $\frak g$ has the $J$ property.
\par
To prove the last statement, we note that,
because $[ J_{\frak g},\frak g_{-2}]=0$, the derivation
$\ad_{\frak g}(J_{\frak s})$ acts on $\frak l_{-2}$
as $-k\cdot \Id_{\frak l_{-2}}$, with $k\in\Bbb K$.
The restriction of 
$\ad_{\frak g}(J_{\frak s})$ to 
$\frak l$ is an endomorphism with zero trace. Since any two 
eigenvalues of $\ad(\Js)$ differ by an integral multiple of 
$i$, all, including $k$,
are purely imaginary. 
\par
In particular, when $\Bbb K=\Bbb R$, we obtain that 
$J_{\frak g}=J_{\frak s}$.
\par
If $\frak l$ is quaternionic, we consider
$Z=\Jg-\Js$. The restriction of
$\ad_{\frak g}(Z)$  to $\frak l$ has purely imaginary
eigenvalues, and
therefore $Z$ belongs to $\frak a\simeq\frak{su}(2)$ (see
Theorem \ref{thm:irredstruct}). Since $\frak a\subset\frak g_0$, we have
$[\Jg,\frak a]=0$ and hence $[Z,\frak a]=-[\Js,\frak a]=0$ shows that $Z=0$.
\end{proof}
\par
\smallskip
Assume that both $\frak s$ and $\frak l$ are complex, i.e. obtained
from a complex Lie algebra and a complex representation by restriction
of the base field to $\Bbb R$. Fix a Cartan subalgebra $\frak h$ of
$\frak s$, with $E_{\frak s}\in\frak h$, so that $\frak h\subset\frak s_0$
and $J_{\frak s}\in\frak h$. We shall denote by $\Cal R$ 
and $\Cal P$ the root
system of $\frak s$ and the set of weights
of $\frak l$ with respect to $\frak h$, respectively.\par
The eigenspaces 
$\frak s^{\alpha}=\{X\in\frak s\mid [H,X]=\alpha(H)X
\;\forall H\in\frak h\}\subset\frak s$,
for $\alpha\in\Cal R$,
and $\frak l^{\lambda}=\{Y\in\frak l\mid [H,Y]=\lambda(H)Y\;
\forall H\in\frak h\}\subset\frak l$, for $\lambda\in\Cal P$,
consist of homogeneous vectors and therefore we have 
partitions:
\[
\Cal R=\bigcup_{p\in\Bbb Z}{\Cal R_p}\quad\text{and}\quad
\Cal P=\bigcup_{p\in\Bbb Z}{\Cal P_p}
\]
where $\Cal R_p$ is the set of the roots $\alpha\in\Cal R$ such that
$\frak s^\alpha\subset\frak s_p$ and $\Cal P_p$ that of the weights
$\lambda\in\Cal P$ for which $\frak l^\lambda\subset\frak l_p$. 
Since $J_{\frak s}\in\frak h$, the set $\Cal R_{-1}$ further
decomposes into
\[
\Cal R_{-1}=\Cal R^{1,0}\cup\Cal R^{0,1}
\]
where $\Cal R^{1,0}=\{\alpha\in\Cal R_{-1}\mid \alpha(J_{\frak s})=i\}$
and 
$\Cal R^{0,1}=\{\alpha\in\Cal R_{-1}\mid \alpha(J_{\frak s})=-i\}$.
We set 
\[
\left\{\begin{aligned}
\frak s^{1,0}&=\{X\in\frak s_{-1}\mid [J_{\frak s},X]=iX\}\, ,\\
\frak s^{0,1}&=\{X\in\frak s_{-1}\mid [J_{\frak s},X]=-iX\}\, ,\\
\frak l^{1,0}&=\{Y\in\frak l_{-1}\mid JY=iY\}\, ,\\
\frak l^{0,1}&=\{Y\in\frak l_{-1}\mid JY=-iY\}\, ,\\
\frak g^{1,0}&=\frak s^{1,0}\oplus\frak l^{1,0}\, ,\\
\frak g^{0,1}&=\frak s^{0,1}\oplus\frak l^{0,1}\, .
\end{aligned}\right.
\]
The integrability condition for the partial complex structure $J$
can be rewritten in the form:
\begin{equation}\label{eq:integrability}\tag{$\diamondsuit$}
[\frak g^{1,0},\frak g^{1,0}]=\{0\}\quad\text{and}\quad
[\frak g^{0,1},\frak g^{0,1}]=\{0\}\, .
\end{equation}
Then we obtain:
\begin{lem2}\label{lem:nostring}
Let $\lambda\in\Cal P_{-1}$.  Then:
\begin{itemize}
\item If $\alpha\in\Cal R_{-1}$ and $\lambda+\alpha\in\Cal P$,
then $\lambda-\alpha\notin\Cal P$;
\item If $\alpha\in\Cal R_{-1}$ and $\lambda-\alpha\in\Cal P$,
then $\lambda+\alpha\notin\Cal P$;
\item $\rm{dim}_{\Bbb C}{\frak l^{\lambda}}=1$ and either
$\frak l^{\lambda}\subset\frak l^{1,0}$, or 
$\frak l^{\lambda}\subset\frak l^{0,1}$.
\end{itemize}
\end{lem2}
\par
\begin{proof}
Let $\alpha\in\Cal R_{-1}$ be fixed, and let $X_{\alpha}$ be a nonzero
element of $\frak s^{\alpha}$. If $X_{\alpha}\in\frak s^{1,0}$ and
$A\in\frak l_0$, then $[X_\alpha,A]\in\frak l^{1,0}$. Indeed we have
$J([X_\alpha,A])=[JX_\alpha,A]=[[J_{\frak s},X_{\alpha}],A]=i[X_\alpha,A]$.
Analogously, if $X_{\alpha}\in\frak s^{0,1}$ and
$A\in\frak l_0$, then $[X_\alpha,A]\in\frak l^{0,1}$.
Therefore $[X_\alpha,X_\alpha,A]=0$ for all $A\in\frak l_0$ 
and $X_\alpha\in\frak s^\alpha$ with $\alpha\in\Cal R_{-1}$. It follows
that $\Cal P$ does not contain any string of the form
$\lambda-\alpha,\lambda,\lambda+\alpha$ with $\alpha\in\Cal R_{-1}$.
\par
On the other hand, since $\frak s\oplus\frak l$ is nondegenerate,
there exists at least one $\alpha\in\Cal R_{-1}$ such that
$\lambda+\alpha\in\Cal P$. Thus $\Cal P$ contains a maximal string of
the form $\lambda,\lambda+\alpha,\hdots,\lambda+p\alpha$, with $p\geq 1$.
It follows that $\frak l^{\lambda}$ is one dimensional and, because
$\frak l^{\lambda}$ is $J$-invariant (indeed $J$ commutes
with the action of $\frak h$ on $\frak l_{-1}$) we conclude, 
from $[X_\alpha,JY]=-[JX_\alpha,Y]=-[[J_\frak s,X_\alpha],Y]=
-\alpha(J_\frak s)[X_\alpha,Y]$, that
$JY=-\alpha(J_\frak s)Y$.
\end{proof}
\par
In particular we obtain a partition
\[
\Cal P_{-1}=\Cal P^{1,0}\cup\Cal P^{0,1}\]
where
\[
\Cal P^{1,0}=\{\lambda\in\Cal P_{-1}\mid
\frak l^{\lambda}\subset\frak l^{1,0}\}\quad\text{and}\quad
\Cal P^{0,1}=\{\lambda\in\Cal P_{-1}\mid
\frak l^{\lambda}\subset\frak l^{0,1}\}\, .
\]
\par
\begin{lem2}\label{lem:liemonom2}
Assume that
$\frak l$ is irreducible and fix $Y\in\frak l_{-2}\setminus\{0\}$.
Then $\frak l_{-2}$ is generated by
the Lie monomials of the form 
\begin{equation}\label{eq:monomials}\tag{$*$}
L=[Z_1,\dots,Z_h,X_1,\Xi_1\dots,X_k,\Xi_k,Y]
\end{equation}
where $Z_j\in\frak s_0$, for $j=1,...,h$, and 
$X_{j}\in\frak s_{-1}$, $\Xi_{j}\in\frak
s_{1}$ for all $j=1,...,k$.
\end{lem2}
\begin{proof}
Since $\frak s$ is transitive,  $\frak s_{-1}\cup\frak s_1$ is a set of
generators of $\frak s$. Therefore, the 
homogeneous Lie monomials $[W_1,\dots,W_\ell,Y]$
with $W_1,\dots,W_\ell\in\frak s_{-1}\cup\frak s_1$ generate $\frak l$. In
particular, those of degree $0$ generate $\frak l_{-2}$. 
By reordering it is
easy to show that these 
are linear combinations of Lie monomials of the form
\eqref{eq:monomials}.
\end{proof}
\begin{lem2}
Let  
$\alpha\in\Cal R^{1,0}$ and $\beta\in\Cal R^{0,1}$,
and 
$X_{\alpha}\in\frak s^{\alpha}$,
$\Xi_{-\alpha }\in\frak s^{-\alpha}$,
$X_{\beta}\in\frak s^{\beta}$,
$\Xi_{-\beta }\in\frak s^{-\beta}$,
$\Xi\in\frak s^{-\beta}$. Then
for all 
$Y\in\frak l_{-2}$, we have
$[X_{\alpha},\Xi_{-\beta},Y]=[X_{\beta},\Xi_{-\alpha},Y]=0$.
\end{lem2}
\begin{proof}
We can assume that $Y\in\frak l^\lambda$ for some weight $\lambda$.
Let us show that $[X_{\alpha},\Xi_{-\beta},Y]=0$.
If $[\Xi_{-\beta},Y]=0$, we have nothing to prove. 
If $[\Xi_{-\beta},Y]\neq 0$, then
$\lambda-\beta$ is a weight of $\frak l$. Since the weights of $\frak l$
contain the string $\lambda$, $\lambda-\beta$, there exist $Y'\in\frak
l^{\lambda-\beta}$ and $X'\in\frak s^{\beta}$ such that
\[
[X',Y]\neq 0\, .
\]
As $\frak s^{\beta}\subset\frak s^{0,1}$, it follows from equation
\eqref{eq:integrability} that $\frak l^{\lambda-\beta}\subset\frak
l^{1,0}$. Hence 
$[\Xi_{-\beta},Y]\in\frak l^{1,0}$ and therefore, using again equation
\eqref{eq:integrability}, we obtain $[X_{\alpha},\Xi_{-\beta},Y]=0$.
\par
The other equality is proved in the same way.
\end{proof}
\par
\smallskip
\begin{proof}[Proof of Theorem \ref{teo:Jprop}] 
\emph{The complex case.}
We can assume that $\frak g$ is proper. Suppose first that 
$\frak l$ is irreducible, and let
$Y\neq 0$ be an eigenvector of $\Js$ in $\frak l_{-2}$, with 
$[\Js,Y]=kY$. By
Lemma \ref{lem:liemonom2}, $\frak l_{-2}$ is generated by the Lie monomials:
\[
L=[Z_1,\dots,Z_h,X_1,\Xi_1\dots,X_k,\Xi_k,Y]
\]
where $Z_1,\dots,Z_h\in\frak s_0$ and moreover
\[\begin{aligned}
X_i\in\frak s^{\alpha_i}&\quad\text{for a root}&\alpha_i\in\Cal R_{-1}\, ,\\
\Xi_i\in\frak s^{-\beta_i}&\quad\text{for a root}&\beta_i\in\Cal R_{-1}\, .
\end{aligned}\]
By the preceding lemma, $L=0$ unless
\[
\alpha_i(\Js)-\beta_i(\Js)=0\quad\text{for every}\quad i=1,\dots,k\, .
\]
Then, if $L\neq 0$ we obtain $[\Js,L]=kL$.
This shows that $\ad(\Js)|_{\frak l_{-2}}$ is a multiple of the identity and
therefore, by Lemma \ref{lem:multidentity}, 
$\frak g$ has the $J$ property.
\par
If $\frak l$ is reducible, $\frak l=\bigoplus_{i=1}^n\frak l_i$, we can
apply the argument 
above to each $\frak s\oplus\frak l_i$: we find that
$J_i=\Js-k_i\pi_i$ for each $i$, where $\pi_i$ is the 
canonical projection
onto $\frak l_i$. Then $J_{\frak g}=\Js-\sum_{i=1}^nk_i\pi_i\in\frak g$ 
defines the partial complex structure of $\frak s\oplus\frak l$.
\par\smallskip
\emph{The general case.}
If $\frak l$ contains irreducible $\frak s$-submodules that are not of
the complex type, we apply the argument above to the complexification
$\Hat{\frak g}=\frak g(\Hat{\frak s}\oplus\Hat{\frak l})$.
In particular we obtain that $\ad(\Js)$ is a multiple of the identity
on the subspace of degree $-2$ of each graded irreducible component of $\Hat{\frak l}$.
If $\frak l^{(i)}$ is a graded irreducible component of $\frak l$
that is of the real type, then its
complexification $\Hat{\frak l}^{(i)}$ is an irreducible
component of $\Hat{\frak l}$. Then by the 
invariance under conjugation it
follows that $\ad(\Js)$ is actually zero on $\frak l^{(i)}_{-2}$.
\par
In the same way, if
$\frak l^{(i)}$ is a graded irreducible component of $\frak l$
that is of the quaternionic type, the complexification 
$\Hat{\frak l}^{(i)}$ is the direct sum of two 
isomorphic simple complex $\Hat{\frak s}$-modules 
$\frak v$ and $\bar{\frak v}$, that are exchanged by the
conjugation. Since $\ad(\Js)$ defines 
on $\frak v_{-2}$ and $\bar{\frak v}_{-2}$ the same purely
imaginary multiple of the identity, this must coincide 
with its conjugate. Thus it is zero and the partial complex structure
$J$ on $\frak l^{(i)}_{-1}$ is defined by
the restriction to $\frak l^{(i)}_{-1}$ of $\ad(\Js)$.\par
This shows that, for the element $J_{\frak g}$, obtained for the 
complexification 
$\Hat{\frak g}=\frak g(\Hat{\frak s}\oplus\Hat{\frak l})$
using the argument in the first part of the proof, we have
$\ad(J_{\frak g})=\ad(\Js)$ on the complexifications
$\Hat{\frak l}^{(i)}$ of all 
graded irreducible components ${\frak l}^{(i)}$ that are of
the real or of the quaternionic type. Therefore
$J_{\frak g}$ belongs to $\frak g(\frak s\oplus\frak l)$.
This completes the proof.
\end{proof}
\par
\begin{ex2}
Let $\frak s=\frak{sl}(2,\Bbb C)$, with the same gradation and $CR$ 
structure as in Example \ref{ex:sl2c1}.  Let $n\geq 2$ and let
$\frak l^{n}$ be the 
irreducible representation of $\frak s$ of dimension $n$.  Lemma 
\ref{lem:nostring} implies that there exists a unique gradation and 
$CR$ structure on $\frak l^{n}$ that makes it into a nondegenerate 
graded $CR$ representation of $\frak s$.  It is obtained by assigning 
degree $-1$ to the minimal weight vector, so that $\frak l^{n}=\frak
l^{n}_{-}$. 
\par
The algebra $\frak s\oplus\frak l^{n}$ is a partial Levi--Tanaka 
extension if $n\geq 2$.  It is semisimple, with $\frak 
g(\frak s\oplus\frak l^{2})\simeq\frak{sl}(3,\Bbb C)$, if $n=2$, and it is
proper, with $\frak g(\frak s\oplus\frak l^{n})=\frak s\oplus\frak 
l^{n}\oplus\Bbb C\pi$, if $n \geq 3$.
\end{ex2}

\subsection{$CR$ admissibility of a graded representation}
\par
Let $\frak s$ be a semisimple graded $CR$ algebra. In this section
we investigate the conditions for a graded $CR$ module $\frak l$
to be a $CR$ representation of $\frak s$.\par
Recall that, if $\frak s$ is Levi--Tanaka, all Levi--Tanaka extensions
of $\frak s$ are proper.
 \par
\begin{lem2}\label{lem:nondeg}
Let $\frak s$ be a fundamental graded Lie algebra and 
$\frak l=\oplus_{p\geq -\nu}{\frak l_p}$, 
with $\frak l_{-\nu}\neq\{0\}$, a 
graded
irreducible representation of $\frak s$. Then for every  
$Y\in\frak l\setminus\frak l_{-\nu}$
there exists $X\in\frak s_{-1}$ 
such that
$[X,Y]\neq 0$.
\end{lem2}
\begin{proof}
It suffices to consider a
$Y\in\frak l_p\setminus\{0\}$, with $p>-\nu$. 
From Lemma \ref{lem:liemonom} we know that 
$\frak l$ is
generated by decreasing homogeneous Lie monomials in $Y$. 
Since $\frak l_{-\nu}\neq\{0\}$,
there are homogeneous elements $Z_1,\hdots,Z_k$ 
of $\frak s_-$ such that $0\neq
[Z_1,\hdots,Z_k,Y]\in\frak l_{-\nu}$. 
In particular $[Z_k,Y]\neq 0$ and,
because $\frak s_-$ is generated by $\frak 
s_{-1}$, 
there is $X$ in $\frak s_{-1}$ such that $[X,Y]\neq 0$.
\end{proof}
\par
We have the following:
\par
\begin{prop2}\label{prop:char1}
A partial Levi--Tanaka extension of a semisimple graded $CR$
algebra $\frak s$ is completely determined by 
the data of a representation $\frak l$,
of a decomposition $\frak l=\oplus{\frak l^{(i)}}$ of
$\frak l$ into a direct sum of irreducible $\frak s$ modules, and,
for each irreducible component $\frak l^{(i)}$,
of its kind $\mu_i=\mu(\frak l^{(i)})$.
\par
Vice versa, given a semisimple Levi--Tanaka algebra $\frak s$ and a 
graded
representation $\frak l$, the extension of $\frak s$ by $\frak l$ is a 
partial
Levi--Tanaka extension if and only if every irreducible 
graded component 
$\frak l^{(i)}$ of $\frak l$ satisfies the following conditions:
\begin{itemize}
\item [($i$)] $\frak l^{(i)}_{-1}\neq 0$, i.e. $\nu(\frak l^{(i)})\geq -1$; 
\item [($ii$)] $\frak l^{(i)}_{-2}\neq 0$, i.e. $\mu(\frak l^{(i)})\geq 2$;
\item [($iii$)] $\frak l^{(i)}$ is $CR$.
\end{itemize}
\end{prop2}
\par
\begin{proof}
The first statement is a direct consequence of Lemma \ref{lem:grad} 
and Lemma
\ref{lem:unique}.
\par
Conditions ($i$), ($ii$) 
and ($iii$) in the second statement are clearly necessary, and 
the previous lemma tells us that, when conditions ($i$) and ($ii$) are
satisfied, $\frak{s\oplus l}$ is transitive and nondegenerate. 
\par
We only need to check that in this case $\frak l$ is also fundamental,
and 
to this aim 
it suffices to consider the case where $\frak l$ is irreducible.\par
Suppose that $\frak l$ is irreducible and
fix $Y\in\frak l_{-1}$. By Lemma \ref{lem:liemonom}, 
the subspace $\oplus_{p<-1}{\frak l_{p}}$ is generated by
homogeneous Lie monomials in $Y$.
Since $\frak s$ is generated by its elements of degree $0$, $1$ and
$-1$, we can use as generators homogeneous 
Lie monomials that only contain elements of
these degrees.\par
Rearranging their terms, 
we obtain a set of generators of
$\oplus_{p<-1}{\frak l_{p}}$ 
that are homogeneous Lie monomials of the form:
\[
[X_1,\dots,X_k,\Xi_1,\dots,\Xi_l,A_1,\dots,A_h,Y]
\]
where $X_i\in\frak s_{-1}$, $\Xi_i\in\frak s_1$ and $A_i\in\frak s_0$ 
for
all $i$, and
$k>l+1$. Every such monomial can be written as 
\[
[X_1,\dots,X_{k-l},W]
\]
where
\[
W=[X_{k-l+1},\dots,X_k,\Xi_1,\dots,\Xi_l,A_1,\dots,A_h,Y]\in\frak l_{-1}\, .
\qedhere\]
\end{proof}
\par
By the previous proposition, to characterize the
homogeneous $CR$ representations $\frak l$ of $\frak s$, 
it will suffice to
restrict our consideration to the case of an irreducible $\frak l$.
\par
\begin{prop2}\label{prop:char2}
Let $\frak s$ be a semisimple graded $CR$ algebra and $\frak l$ a
graded irreducible representation of $\frak s$. Then $\frak{s\oplus l}$ 
is a partial Levi--Tanaka extension of $\frak s$ if and only if
\begin{itemize}
\item [($i$)] $\frak l_{-1}\neq 0$; 
\item [($ii$)] $\frak l_{-2}\neq 0$;
\item [($iii$)] there exists $k\in\Bbb R$ such that $\ad(\Js-ik\pi)=0$ on 
$\frak l_0\oplus \frak l_{-2}$. 
\end{itemize}
\end{prop2}
\par
\begin{proof}
The conditions are necessary in view of Proposition \ref{prop:char1},
Theorem \ref{teo:Jprop} and  Lemma \ref{lem:multidentity}, since by
Lemma 3.11 in
\cite{MN97}, the element $J_{\frak g}$ of the Levi--Tanaka prolongation
of $\frak s\oplus\frak l$ that defines the partial complex structure
satisfies $[J_{\frak g},\frak l_0\oplus\frak l_{-2}]=\{0\}$.
\par
Assume vice versa that conditions ($i$),($ii$) and  ($iii$) are
satisfied and define (we keep for $\pi$ and $\frak a$ the notation
of the preceding sections)
$\Tilde J=J_{\frak s}-ik\pi\in\frak s\oplus\frak l\oplus\frak a$. 
Then
\[
[\Tilde J,X,Y]=[[\Tilde J,X],Y]+[X,[\Tilde J,Y]]=
\begin{cases}
0               &\text{if}\ X\in\frak s_{-1},\, Y\in\frak l_{-1}; \\
[X,[\Tilde J,Y]]&\text{if}\ X\in\frak s_{0},\, Y\in\frak l_{-1}; \\
[JX,Y]          &\text{if}\ X\in\frak s_{-1},\, Y\in\frak l_{0}; 
\end{cases}
\]
and, for $X\in\frak s_{-1}$, $Y\in\frak l_{-1}$,
\begin{multline*}
[X,[\Tilde J,[\Tilde J,Y]]]=[\Tilde J,[X,[\Tilde J,Y]]]-[JX,[\Tilde 
J,Y]= \\
=0-[\Tilde J,[JX,Y]]+[JJX,Y]=-[X,Y].
\end{multline*}
Thus $\ad(\Tilde J)$ defines a partial complex structure on 
$\frak{s\oplus l}$ and, by Proposition \ref{prop:char1}, 
$\frak{s\oplus l}$ is a partial Levi--Tanaka extension of $\frak s$.
\end{proof}
\par
Finally we give a characterization of the graded $CR$ representations
$\frak l$ of a semisimple graded $CR$ algebra $\frak s$ in terms
of the weights of the (complexification of the) $\frak s$-module $\frak l$.
\par
We recall from \cite{MN97} that the complexification of a Levi--Tanaka
algebra is still a Levi--Tanaka algebra and that vice versa a
Levi--Tanaka algebra structure on the complexification $\hat{\frak g}$
of a graded Lie algebra 
$\frak g=\oplus_{p\in\Bbb Z}{\frak g_p}$ 
induces a Levi--Tanaka structure
on $\frak g$ if and only if $\frak g_{-1}$ is invariant for the
complex structure $J$ of $\hat{\frak g}_{-1}$.
\par
We use the same notation of Section 
\ref{sect:CRstruct}. Moreover, we denote by $\Cal B$ a system of
simple vectors of $\Cal R$, such that the 
eigenspaces of the negative roots are
contained in $\Hat{\frak s}_+$ and we have a partition
$\Cal B=\Cal B_0\cup\Cal B^{1,0}\cup\Cal B^{0,1}$,
with $\Cal B_0=\Cal B\cap\Cal R_0$, 
$\Cal B^{1,0}=\Cal B\cap\Cal R^{1,0}$,
$\Cal B^{0,1}=\Cal B\cap\Cal R^{0,1}$.

\par
\begin{thm2}\label{thm:admissibility}
Let $\frak s$ be a semisimple graded $CR$ algebra 
of the complex type and $\frak l$ 
a complex irreducible graded representation of $\frak s$, such 
that $\frak l_{-1}\neq 0$ and $\frak l_{-2}\neq 0$. Then the 
following conditions are equivalent:
\begin{itemize}
\item[($i$)] $\frak s\oplus\frak l$ is a partial Levi--Tanaka extension of 
$\frak s$;
\item[($ii$)] $\ad(J_{\frak s})$ acts as a scalar multiple of the 
identity on
$\frak l_0\oplus\frak l_{-2}$;
\item[($iii$)] if $\lambda\in\Cal P_{-2}$, $\alpha,\alpha'\in\Cal R^{1,0}$, 
$\beta,\beta'\in\Cal R^{0,1}$, then none of the weights
$\lambda-\alpha-\alpha'$, $\lambda-\alpha+\beta$, $\lambda-\beta-\beta'$,
$\lambda+\alpha-\beta$ belongs to $\Cal P$; 
\item[($iv$)] There exists a partition $\Cal P_{-1}=\Cal 
P^{1,0}\cup\Cal P^{0,1}$ such that
\begin{itemize}
\item[(\textup{a})] $(\Cal P^{1,0}\pm\Cal B_0)\cap\Cal P\subset\Cal P^{1,0}$ 
and 
$(\Cal P^{0,1}\pm\Cal B_0)\cap\Cal P\subset\Cal P^{0,1}$,
\item[(\textup{b})] $(\Cal P^{1,0}+\Cal B^{1,0})\cap\Cal P=\emptyset$, 
$(\Cal P^{1,0}-\Cal B^{0,1})\cap\Cal P=\emptyset$, $(\Cal 
P^{0,1}+\Cal B^{0,1})\cap\Cal P=\emptyset$, $(\Cal P^{0,1}-\Cal 
B^{1,0})\cap\Cal P=\emptyset$.
\end{itemize}
\end{itemize}
\end{thm2} 

\begin{proof}
($i$) $\Leftrightarrow$ ($ii$).\quad
 This is Proposition \ref{prop:char2} in the 
complex case.
\par\noindent
($ii$) $\Rightarrow$ ($iii$).\quad In fact we obtain
$(\lambda-\alpha-\alpha')(J_{\frak s})=\lambda(J_{\frak s})-2i$,
$(\lambda-\beta-\beta')(J_{\frak s})=\lambda(J_{\frak s})+2i$,
$(\lambda-\alpha+\beta)(J_{\frak s})=\lambda(J_{\frak s})-2i$,
$(\lambda+\alpha-\beta)(J_{\frak s})=\lambda(J_{\frak s})+2i$.
If any of these weights belongs to $\Cal P$, statement ($ii$)
is contradicted.
\par\noindent
($iii$) $\Rightarrow$ ($iv$). \quad 
For each $\lambda\in\Cal P_{-2}$ define 
$\Cal P_{\lambda}^{0,1}=\left(\lambda - \Cal R^{1,0}\right)\cap\Cal P$ 
and 
$\Cal P_{\lambda}^{1,0}=\left(\lambda - \Cal R^{0,1}\right)\cap\Cal P$. 
Let 
$\Cal P^{0,1}=\bigcup_{\lambda}\Cal P_{\lambda}^{0,1}$ and
$\Cal P^{1,0}=\bigcup_{\lambda}\Cal P_{\lambda}^{1,0}$. This is a partition 
of $\Cal P_{-1}$. Indeed every weight in $\Cal P_{-1}$ is contained 
either
in $\Cal P^{0,1}$ or in $\Cal P^{1,0}$, because $\frak l$ is 
nondegenerate. The fact that $\Cal P^{0,1}\cap\Cal P^{1,0}=\emptyset$ 
and conditions ($iv$,{a}) and ($iv$, b) follow  from condition ($iii$).
\par\noindent
($iv$) $\Rightarrow$ ($i$).\quad
 Define a complex structure $J$ on $\frak l_{-1}$ 
by imposing 
that $J$ has eigenvalues $i$ and $-i$ on weight spaces corresponding 
to weights in 
$\Cal P^{1,0}$ and $\Cal P^{0,1}$, respectively. This 
defines a complex structure that
is compatible with the $CR$ structure of $\frak s$ and thus defines a 
$CR$ structure on $\frak l$. Proposition \ref{prop:char1} then implies that 
$\frak s\oplus\frak l$ is a Levi--Tanaka extension of $\frak s$.
\end{proof}

\section{Examples}
\par
In this section we give some examples
of the Levi--Tanaka extensions of semisimple graded $CR$ 
algebras.
\par
\subsection{The Standard manifold of a Levi--Tanaka extension}
\par
To every Levi--Tanaka algebra $\frak g$ corresponds 
a homogeneous 
$CR$ manifold $M(\frak g)$, called the \emph{standard $CR$ manifold} 
associated to $\frak g$ (see \cite{MN2000} for more details). To 
construct $M(\frak g)$, consider the simply connected Lie group $G$ 
with Lie algebra $\frak g$ and let $G_+$ be the connected analytic 
subgroup of $G$ whose Lie algebra is $\frak g_+$. The analytic 
subgroup $G_+$ is closed, and $M(\frak g)=G/G_+$ is a smooth 
manifold. The tangent space at $o=eG_+$ is isomorphic to $\frak 
g/\frak g_+$ and can be identified with $\frak g_-$. In this way a 
partial complex structure is defined on $T_oM$; via left 
translations by the elements of $G$ one defines on $M(\frak g)$ a 
$G$-invariant $CR$ structure. 
In this way we
obtain a complete, Levi-nondegenerate $CR$ manifold, for which
$G$ 
is exactly 
the identity component of 
the group of $CR$ automorphisms. \par
Furthermore $M(\frak g)$ is compact if and only if $\frak g$ is 
semisimple.
\par
In this context Levi--Tanaka extensions correspond to homogeneous 
vector bundles: denote by $S$ and $S_+$ the connected analytic 
subgroups of $G$ corresponding to $\frak s$ and $\frak s_+$. In 
\cite{MN2000} it is proved that $G/G_+$ is a Mostow fibration on 
$S/S_+$. The fiber over $o=eS_+$ is naturally identified with 
$\frak{l/l_+}$. Being $\frak l$ a linear representation of $\frak s$, 
and consequently of $S$, $G/G_+$ has a canonical structure of an 
$S$-homogeneous vector bundle over $S/S_+$ and the isomorphism
\[
G/G_+ \simeq S \times_{S_+} \left(\frak l/\frak l_+\right)
\]
is $S$-equivariant. The isomorphism $\frak l/\frak l_+\simeq\frak 
l_-$ induces a partial complex structure  on every fiber that 
coincides with the $CR$ structure induced by the inclusion in $G/G_+$.
\par\medskip
\subsection{The adjoint representation I}
\par
If $\frak l$ is the adjoint representation of $\frak s$ then 
$\frak{s\oplus l}$ is always a partial Levi--Tanaka extension of 
$\frak s$ (the proof is straightforward) and the resulting vector 
bundle is isomorphic to the tangent bundle of $S/S_+=
M(\frak s)$, with the natural $CR$ structure.
\par
\subsection{The adjoint representation II}
If all simple ideals of 
$\frak s$ are fundamental graded transitive $CR$ algebras (not 
necessarily Levi--Tanaka), then the adjoint action of $\Js$ 
defines a complex structure on $\frak s_{1}$, because the eigenvalues 
of $\Js$ on $\frak s_{1}$ can only be $\pm i$.  \par If $X\in\frak 
s_{2}$ and $Y\in\frak s_{-2}$ then
\[
\kappa([\Js,X],Y)=-\kappa(X,[\Js,Y])=0
\]
where $\kappa$ denotes the Killing form, and hence $[\Js,\frak s_{2}]=0$.  
Furthermore, $[\Js,\frak s_{0}]=0$ because the action of $\frak s_{0}$ 
is compatible with the $CR$ structure.  
\par
Let $\frak l$ be the adjoint representation, with a ``shifted'' 
gradation: $\frak l_{d}\simeq\frak s_{d+2}$.  The relations above 
imply that $\ad(\Js)$ defines a $CR$ structure on $\frak l$, that is 
compatible with the action of $\frak s$.  It is easy  to check 
that the resulting algebra is transitive, fundamental and nondegenerate.
Therefore $\frak s\oplus\frak l$ is a partial Levi--Tanaka extension of 
$\frak s$.
\par
\subsection{Levi--Tanaka extensions of $\frak{su}(1,2)$} 
\par
In this subsection
we classify all Levi--Tanaka extensions of $\frak{su}(1,2)$. 
To this aim, passing to complexifications, we first classify all 
complex 
Levi--Tanaka extensions of $\Hat{\frak s}=\frak{sl}(3,\mathbb C)$. 
\par
We know from \cite{MNssempl} that $\Hat{\frak s}$ admits an 
essentially unique structure of Levi--Tanaka algebra. More precisely,
the gradation and the complex structure are defined by the
inner derivations associated to the matrices: 
\[
E_{\Hat{\frak s}}=
\begin{pmatrix}
1& 0 & 0 \\
0 & 0 & 0 \\
0 & 0 & -1\\
\end{pmatrix}\qquad\text{and}\qquad
J_{\Hat{\frak s}}=
\begin{pmatrix}
\frac{1}{3}i & 0 & 0 \\
0 & -\frac{2}{3}i & 0 \\
0 & 0 & \frac{1}{3}i \\
\end{pmatrix}.
\]
\par
\begin{figure}[ht]
        \centering
        \includegraphics{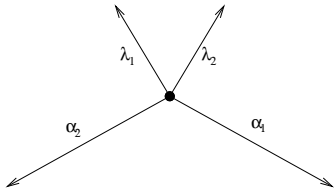}
        \caption{Roots and weights of $\frak{sl}(3,\mathbb C)$}
        \label{fig:basis}
\end{figure}

Fix the Cartan subalgebra $\Hat{\frak h}$ of  diagonal matrices. 
Let $\Cal B=\left\{\alpha_1,\alpha_2\right\}$ be the standard basis 
of the root 
system $\Cal R$ of $\Hat{\frak s}$. Roots and weights of $\Hat{\frak s}$ 
are all elements of a two dimensional Euclidean space, that we 
identify with the standard Euclidean space $\mathbb R^2$. If 
$\omega_1$ and $\omega_2$ are the fundamental dominant weights, 
every irreducible representation 
$\Hat{\frak l}$ is characterized by its maximal weight $\omega$, 
that can be written as $\omega=k_1\omega_1+ k_2\omega_2$, with $k_1$, $k_2$
non negative integers; the corresponding representation $\Hat{\frak l}$ is denoted by 
$\Gamma_{k_1,k_2}$.
Figure \ref{fig:basis} shows the basis $\Cal B$ and the 
fundamental dominant weights.
\par
\par
\begin{figure}[ht]
        \centering
        \includegraphics{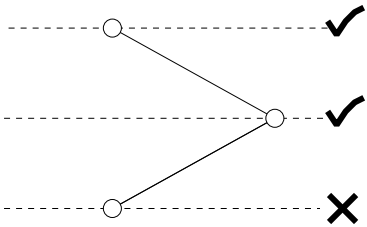}  
        \caption{Weights of $\Gamma_{1,0}$}
        \label{fig:10}
\end{figure}

The representation $\Gamma_{1,0}$ is the standard representation of 
$\frak{sl}(3,\mathbb C)$. Its weight diagram is depicted in Figure 
\ref{fig:10}
\par
\par
In our pictures weights of the same degree will lay on a same dotted line; we
put a mark on the right of a dotted line if there is a partial Levi--Tanaka
extension by $\Hat{\frak l}$ in which that line corresponds to degree $-1$.
\par
\begin{figure}[ht]
        \centering
        \includegraphics{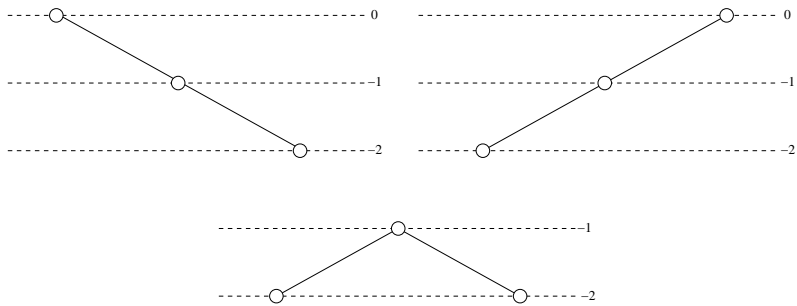} 
        \caption{``Forbidden'' configurations}
        \label{fig:forbidden}
\end{figure}

For $\frak{sl}(3,\mathbb C)$ the equivalent conditions of 
Theorem \ref{thm:admissibility} can be graphically described by 
saying that none of the configurations of Figure \ref{fig:forbidden} could 
appear in the weight diagram of $\Hat{\frak l}$.
\par
\par
\begin{figure}[ht]
        \centering
        \includegraphics{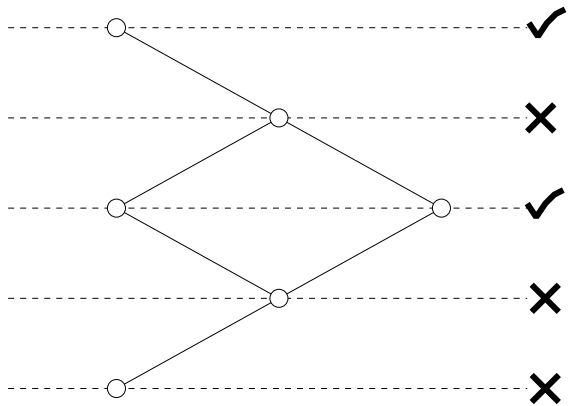}   
        \caption{Weights of $\Gamma_{2,0}$}
        \label{fig:20}
\end{figure}
\begin{figure}[ht]
        \centering
        \includegraphics{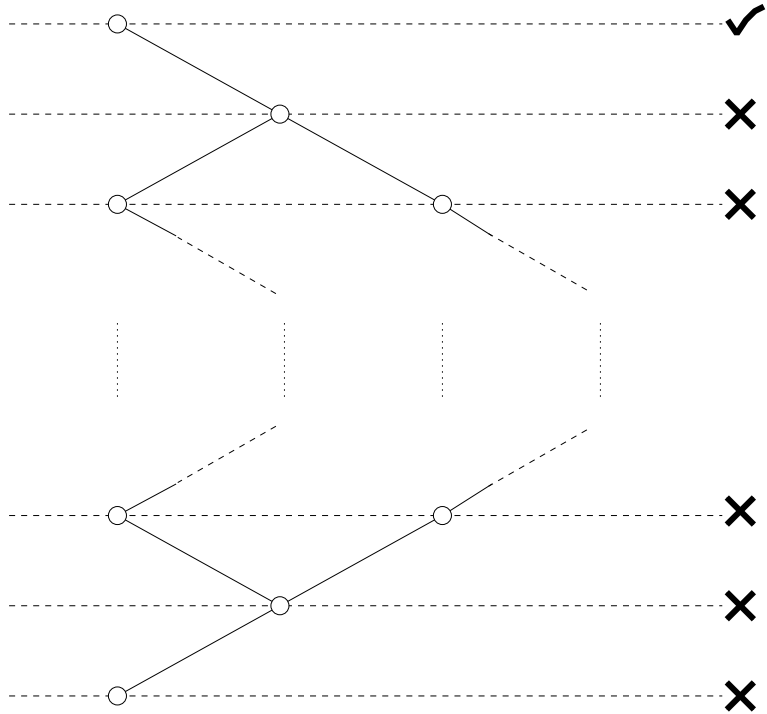}   
        \caption{Weights of $\Gamma_{n,0}$}
        \label{fig:n0}
\end{figure}

The situation for $\Gamma_{2,0}$ is depicted in 
Figure \ref{fig:20}. For representations of type $\Gamma_{n,0}$ with 
$n\geq 3$ it is possible to satisfy the conditions of
Theorem \ref{thm:admissibility} 
only by assigning degree $-1$ to the homogeneous subspace of maximal 
degree (see Figure \ref{fig:n0}).
\par
\par
\par
Representations of type $\Gamma_{0,n}$ are  analogous to those of type 
$\Gamma_{n,0}$.
\par
\begin{figure}[ht]
        \centering
        \includegraphics{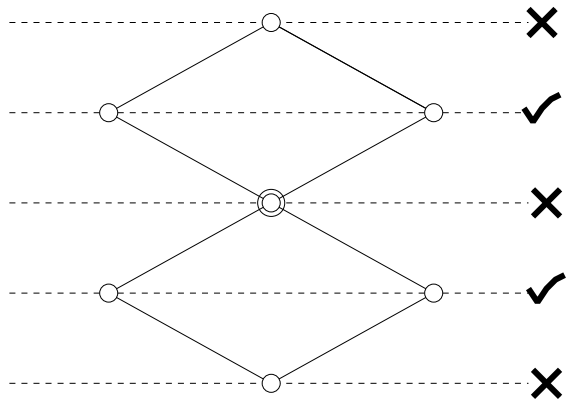}   
        \caption{Weights of $\Gamma_{1,1}$}
        \label{fig:11}
\end{figure}

The only representation $\Gamma_{k_1,k_2}$ with $k_1$ and $k_2$ both positive
that admits a $CR$ structure is the adjoint representation $\Gamma_{1,1}$.
In this case there are exactly two distinct $CR$ structures , that are those
described in the previous subsections and are depicted in Figure \ref{fig:11}.
\par
\par
Now we apply these results to the case of the real algebra 
$\frak s=\frak{su}(1,2)$.
To this aim, we need to consider the $CR$ admissible complex representations
of $\Hat{\frak s}$ as real representations of $\frak s$.
While  $\Gamma_{1,1}$ splits 
into two copies of the adjoint representation $\frak l_{\ad}$, 
the representations $\Gamma_{n,0}$ and $\Gamma_{0,n}$
remain irreducible. Hence the irreducible partial Levi--Tanaka extensions 
of $\frak s$ are $\Gamma_{n,0}$, $\Gamma_{0,n}$ and the adjoint 
representation, with the gradations described above and the corresponding 
unique compatible $CR$ structures.
\par
The representations $\Gamma_{n,0}$ and $\Gamma_{0,n}$ are of the 
complex type, the adjoint representation $\frak l_{\ad}$ is of the 
real type. Their corresponding Levi--Tanaka extensions 
are:
\begin{eqnarray*}
\frak s\oplus\Gamma_{n,0}\oplus\mathbb C\pi, \\
\frak s\oplus\Gamma_{0,n}\oplus\mathbb C\pi, \\
\frak s\oplus\frak l_{\ad}\oplus\mathbb R\pi.
\end{eqnarray*}


\providecommand{\bysame}{\leavevmode\hbox to3em{\hrulefill}\thinspace}

\end{document}